\documentclass{amsart}

\usepackage{amsmath,amssymb,epic}

\newtheorem{theorem}{Theorem}[section]
\newtheorem{proposition}[theorem]{Proposition}
\newtheorem{lemma}[theorem]{Lemma}
\newtheorem{corollary}[theorem]{Corollary}

\theoremstyle{definition}

\theoremstyle{remark}
\newtheorem{remark}[theorem]{Remark}

\numberwithin{equation}{section}

\def\DJ{{\hbox{D\kern-.8em\raise.15ex\hbox{--}\kern.35em}}}
\def\DJo{$\;$\kern-.4em
    \hbox{D\kern-.8em\raise.15ex\hbox{--}\kern.35em okovi\'c}}

\def\al{{\alpha}}
\def\be{{\beta}}

\def\la{{\lambda}}

\def\bZ{{\mbox{\bf Z}}}

\def\bp{{\mbox{\bf p}}}
\def\pB{{\mathcal B}}
\def\pN{{\mathcal N}}
\def\pS{{\mathcal S}}
\def\pT{{\mathcal T}}
\def\pV{{\mathcal V}}

\def\tr{{\rm tr\;}}

\def\Sk{{\mbox{\rm Skew}}}
\def\GL{{\mbox{\rm GL}}}
\def\SL{{\mbox{\rm SL}}}

\def\Ort{{\mbox{\rm O}}}

\renewcommand{\subjclassname}{\textup{2000} Mathematics Subject
Classification}

\begin{document}

\title[Normal forms of skew-symmetric matrices]
{Normal forms for orthogonal similarity classes of
skew-symmetric matrices}

\author[D.\v{Z}. \DJ okovi\'{c}, K. Rietsch and K. Zhao]
{Dragomir \v{Z}. \DJ okovi\'{c}, Konstanze Rietsch and Kaiming Zhao}

\address{Department of Pure Mathematics, University of Waterloo,
Waterloo, Ontario, N2L 3G1, Canada}

\address{Department of Mathematics, King's College London, Strand,
London WC2R 2LS}

\address{Department of Mathematics, Wilfrid Laurier University,
Waterloo, Ontario, N2L 3C5, Canada and Institute of Mathematics,
Academy of Mathematics and System Sciences, Chinese Academy of
Sciences, Beijing 100080, P.R. China}

\email{djokovic@uwaterloo.ca} \email{konstanze.rietsch@kcl.ac.uk}
\email{kzhao@wlu.ca}

\thanks{
The first author was supported by the NSERC Grant A-5285, the
second by a Royal Society Dorothy Hodgkin Research Fellowship, and
the third by the NSERC and the NSF  of China (Grants 10371120 and
10431040). The second author is currently funded by EPSRC grant 
EP/D071305/1.}


\date{}

\begin{abstract}
Let $F$ be an algebraically closed field of characteristic
different from $2$. Define the orthogonal group, $\Ort_n(F)$, as
the group of $n$ by $n$ matrices $X$ over $F$ such that $XX'=I_n$,
where $X'$ is the transpose of $X$ and $I_n$ the identity matrix.
We show that every nonsingular $n$ by $n$ skew-symmetric matrix
over $F$ is
orthogonally similar to a bidiagonal skew-symmetric matrix. In
the singular case one has to allow some 4-diagonal blocks as well.

If further the characteristic is 0, we construct the
normal form for the $\Ort_n(F)$-similarity classes of skew-symmetric
matrices. In this case, the known normal forms
(as presented in the well known book by Gantmacher) are quite
different.

Finally we study some related varieties of matrices. We prove that the
variety of normalized nilpotent $n$ by $n$
bidiagonal matrices for $n=2s+1$ is irreducible of dimension $s$. As
a consequence the skew-symmetric nilpotent $n$ by $n$ bidiagonal matrices
are shown to form a variety of pure dimension $s$.
\end{abstract}
\maketitle \subjclassname{ 15A21, 20G20, 11E57 }

\section{Introduction}

In this note $F$ denotes an algebraically closed field  of
characteristic not $2$. By $M_n(F)$ we denote the algebra of $n$
by $n$ matrices over $F$, by $\GL_n(F)$ the group of invertible
elements of $M_n(F)$, and by $I_n$ the identity matrix of size
$n$. For any matrix $X$, let $X'$ denote the transpose of $X$. The
orthogonal group, $\Ort_n(F)$, is defined as the subgroup of
$\GL_n(F)$ consisting of matrices $X$ such that $XX'=I_n$.

This paper is a sequel to \cite{DZ} where the first and third author constructed 
tridiagonal
normal forms for symmetric matrices under the action of $\Ort_n(F)$.
Here we continue this work and
study the similarity action of $\Ort_n(F)$ on the space
$\Sk_n(F)$ of all $n$ by $n$ skew-symmetric matrices over $F$. We show
that each $A\in\Sk_n(F)$ is orthogonally similar to
the direct sum of blocks $B$ which are either bidiagonal or 4-diagonal.
The latter is needed only if $A$ has at least one pair of
nilpotent Jordan blocks of even size.
Each of the blocks $B$ is similar to either a pair of
Jordan blocks having the same size, $s$, and having eigenvalues
$\la$ and $-\la$, with $s$ even if $\la=0$, or a
single nilpotent Jordan block of odd size.

If $F$ has characteristic $0$, we are able to choose the concrete
representatives, i.e., normal forms for skew-symmetric
matrices, see Theorem \ref{canon}.
In this case, the known normal forms for the
$\Ort_n(F)$-similarity classes of nonsingular skew-symmetric matrices,
as presented in the well known book by Gantmacher \cite{FG},
are not bidiagonal. Hence our new normal
forms are simpler and may be better suited for some applications.
At the end of Section \ref{NF} we pose three open problems.
The proof of Lemma \ref{druga} uses some ideas of Givental and, due
to its length, is given separately in Section 4.

In the last section we study the variety $\mathcal B_n\cap\mathcal N_n$ of nilpotent
$n$ by $n$ bidiagonal matrices with $1$'s along the lower diagonal (normalized
nilpotent bidiagonal matrices), where $n=2s+1$.
We compare $\mathcal B_n\cap \mathcal N_n$ with the variety of
normalized nilpotent tridiagonal  matrices $\mathcal T_n\cap \mathcal N_n$
studied by Kostant in relation
with the Toda lattice \cite{Kos:Toda}, and
we show that $\mathcal B_n\cap \mathcal N_n$ is
irreducible of dimension $s$. These results also have consequences for
the related variety of skew-symmetric nilpotent bidiagonal matrices,
see Corollary~\ref{skewvar}.
Finally we note that the
coordinate ring of $\mathcal B_n\cap\mathcal N_n$
has an interpretation in terms of quantum cohomology of the flag variety
$\SL_n/B$.

The second author thanks the Univeristy of Waterloo for its hospitality 
during the second half of 2005.

\section{Bidiagonal and 4-diagonal representatives}

We first recall the following well known facts (see e.g.
\cite[Theorem 70]{IK} and \cite[Chapter XI, Theorem 7]{FG}).
Gantmacher gives the proof over complex numbers but his argument
is valid for any algebraically closed field of characteristic
different from 2.

\begin{theorem} \label{ort}
If two skew-symmetric matrices $A,B\in M_n(F)$ are similar, then they
are orthogonally similar.
\end{theorem}

Let $t$ be an indeterminate over $F$ and let $i$
denote one of the two square roots of $-1$ in $F$. Thus $i^2=-1$.

\begin{theorem} \label{kose-ort}
A matrix $A\in M_n(F)$ is similar to a skew-symmetric matrix
iff the elementary divisors $(t-\la)^s$ and $(t+\la)^s$
of $A$ come in pairs if $\la\ne0$ or $s$ is even.
\end{theorem}

We are interested in constructing normal forms of skew-symmetric
matrices under orthogonal similarity. The case of complex matrices
is classical and is described in Gantmacher's book \cite{FG}.
Our objective is to construct a normal form which is
almost bidiagonal, unlike the one given by Gantmacher.

The general form of a skew-symmetric bidiagonal matrix is
\[
S=\left[ \begin{matrix}
0 & a_1 & 0 & 0 & \cdots & 0 & 0 & 0 \\
-a_1 & 0 & a_2 & 0 & & 0 & 0 & 0 \\
0 & -a_2 & 0 & a_3 & & 0 & 0 & 0 \\
0 & 0 & -a_3 & 0 & & 0 & 0 & 0 \\
\vdots &  &  &  & &  &  &  \\
0 & 0 & 0 & 0 & & 0 & a_{n-2} & 0 \\
0 & 0 & 0 & 0 & & -a_{n-2} & 0 & a_{n-1} \\
0 & 0 & 0 & 0 & & 0 & -a_{n-1} & 0
\end{matrix} \right].
\]
It is easy to see that if all $a_k$'s are nonzero, then $S$ is
cyclic (i.e., its minimal and characteristic polynomials
coincide). If all $a_k\ne0$, and $a_{k}=1$ for even $k$'s
then we say that this bidiagonal matrix is {\em special}.

Let $s$ be a positive integer.
We introduce a skew-symmetric nilpotent matrix $Q_{4s}$ of size $4s$
having the elementary divisors $t^{2s},t^{2s}$. The easiest way
to understand the structure of this matrix is to take a look at
one example:
\[
Q_{12}=\left[\begin{array}
{ccccccccccccccc}
\cline{2-3}
\multicolumn{1}{c|}{0}
&-1&
\multicolumn{1}{c|}{i}
   &&&&&&&&&&&0\\ %2
   \cline{2-3}
1&0&
0
&&&&&\\ %3
  \cline{4-5}
-i&0&
\multicolumn{1}{c|}{0}
 &-1&\multicolumn{1}{c|}{i}
&& &&&&\\ %4
  \cline{4-5}
&&
1
 &0&0
&& &&&&\\ %5
  \cline{7-8}&&
-i
&0&
0&
 & \multicolumn{1}{|c}{-1}&\multicolumn{1}{c|}{i}
&& &&\\ %6
   \cline{7-8}
&&&&
1
 &&0
&0&&
-i
&&&\\ %7
  &&&&
-i&&
0&
0
 &&-1
&& \\ %8
\cline{7-8}
 &&&&
 &&\multicolumn{1}{|c}{i}
&\multicolumn{1}{c|}{1}&&
0&0&
-i
 &\\ %9
\cline{7-8}
  &&&&&&&&&0&
0
 &-1
&&
\\ %10
 \cline{10-11}
&&&&&&&&
\multicolumn{1}{c|}{}
 &i&\multicolumn{1}{c|}{1}&
0&0&-i\\ %11
\cline{10-11}
  &&&&&&&&&&
&0&
0 &{-1}
 \\ %12
 \cline{12-13}
  0&&&&&&&&&&
\multicolumn{1}{c|}{}&i&
\multicolumn{1}{c|}{1} &{0}
\\ %13
\cline{12-13}
&&&&&&&&&&&&&
\end{array}\right]
\]

In order to define precisely these matrices, let us write $E_{p,q}$
for a square matrix (of appropriate size, $4s$ in this case) all of
whose entries are 0 except the $(p,q)$-entry which is 1. Then we
have
\[ Q_{4s}=X-X', \]
where $X=X_1+X_2$ and
\begin{eqnarray*}
X_1 &=& \sum_{p=1}^s \left(iE_{2p-1,2p+1}-E_{2p-1,2p} \right), \\
X_2 &=& \sum_{p=s+1}^{2s} \left(iE_{2p,2p-2}+E_{2p,2p-1} \right).
\end{eqnarray*}
In the above example, the nonzero entries of $X$ are enclosed
in small boxes.

We shall now prove the above mentioned properties.

\begin{lemma}
The skew-symmetric matrix $Q_{4s}$ is nilpotent and has
elementary divisors $t^{2s},t^{2s}$.
\end{lemma}

\begin{proof}
It is easy to verify that $XX'=0$, and also that $X_1X_2=X_2X_1=0$
and $X_1^{s+1}=X_2^{s+1}=0$.
The matrix $X_1^s$ resp. $X_2^s$ has exactly two nonzero
entries namely $i^{s+1},i^s$ resp. $i^s,i^{s-1}$ and they are
positioned in the middle of the first resp. last row.
It follows easily that $X^{s+1}=0$, $X^s=X_1^s+X_2^s$, and
${X'}^sX^s=0$. By using these relations and $XX'=0$, we obtain that
\begin{eqnarray*}
Q_{4s}^{2s} &=& (X-X')^{2s} \\
&=& X^{2s}-X'X^{2s-1}+{X'}^2X^{2s-2}-\cdots-{X'}^{2s-1}X+{X'}^{2s} \\
&=& (-1)^s{X'}^sX^s = 0.
\end{eqnarray*}
As $Q_{4s}$ is nilpotent of rank $4s-2$, it follows that it
has exactly two Jordan blocks each of size $2s$.
\end{proof}

\begin{remark} The matrix $Q_{4s}$ in the above lemma cannot be replaced by
a bidiagonal skew-symmetric matrix $S$ (of size $n=4s$). Indeed assume 
that $S$ is the bidiagonal matrix displayed above with superdiagonal
entries $a_1,...,a_{n-1}$. Since $S$ must have rank $n-2$, exactly one
of the $a_k$'s is $0$. As it has elementary divisors $t^{2s}, t^{2s}$
we conclude that $a_{2s}=0$. This contradicts Theorem~2.2 as the
$2s\times 2s$ block in the upper left hand corner has size $2s$ and
only one elementary divisor $t^{2s}$.
\end{remark}

We can now state and prove one of our main results.

\begin{theorem} \label{alm-bid} Every skew-symmetric matrix $A\in M_n(F)$
is orthogonally similar to the direct sum
$A_1\oplus A_2\oplus\cdots\oplus A_m$,
where each block $A_k$ is one of the following:

(a) special bidiagonal of size $2s$ with the elementary
divisors $(t-\la)^s,(t+\la)^s$ and $\la\ne0$;

(b) special bidiagonal of size $2s+1$ with the elementary
divisor $t^{2s+1}$;

(c) $Q_{4s}$, with the elementary divisors $t^{2s},t^{2s}$.

\end{theorem}

\begin{proof}
In view of Theorems \ref{ort} and \ref{kose-ort}, it suffices to construct
$A\in\Sk_n(F)$ having only one or two elementary divisors, as specified
in the three cases of the theorem.

Let us start with the case (a). Then $n=2s$ and the
two elementary divisors are $(t-\la)^s$ and $(t+\la)^s$
with $\la\ne0$.
Let $S=S(x_1,x_2,\ldots,x_s,y)$ be the skew-symmetric bidiagonal
matrix whose entries on the first superdiagonal
are the indeterminates
\[ x_1,y,x_2,y,\ldots,x_{s-1},y,x_s \]
in this order.

By using an obvious permutation
matrix $P$, we can transform $S^2$ to obtain direct sum of two
matrices of size $s$ each:
\[ PS^2P^{-1}=S_1\oplus S_2. \]
The matrix $S_1$ resp. $S_2$ is the submatrix of $S^2$ which lies in
the intersection of rows and columns having odd resp. even indices.
These two matrices are tridiagonal and symmetric.
Explicitly, the entries on the first superdiagonal of $S_1$ are
\[ x_1y,x_2y,\ldots,x_{s-1}y; \]
and its diagonal entries are
\[ -x_1^2, -x_2^2-y^2, \ldots, -x_{s-1}^2-y^2, -x_s^2-y^2. \]
The corresponding entries of $S_2$ are
\[ x_2y,x_3y,\ldots,x_sy; \]
and
\[ -x_1^2-y^2, -x_2^2-y^2, \ldots, -x_{s-1}^2-y^2, -x_s^2. \]

We claim that the matrices $S_1=S_1(x_1,\ldots,x_s,y)$
and $S_2=S_2(x_1,\ldots,x_s,y)$ have the same characteristic
polynomial. Let $X$ be the $s\times s$ matrix whose diagonal
entries are $x_1,x_2,\ldots,x_s$,
those on the first subdiagonal are all equal to $-y$, while all
other entries are 0. Then $S_1=-XX'$ and $S_2=-X'X$ and our claim
follows.

Let us write the characteristic polynomial of the matrix
$\la^2 y^2 I_s-S_1(x_1,\ldots,x_s,y)$ as
\[ p(t)=t^s+c_1t^{s-1}+\cdots+c_{s-1}t+c_s, \]
where
$c_k=c_k(x_1,\ldots,x_s,y)$
is a homogeneous polynomial of degree $2k$ in the indicated variables.

Since there are $s+1$ indeterminates
and $F$ is algebraically closed, by \cite[Theorem 4, Corollary 5,
p.57]{IS} the system of homogeneous equations:
$$
c_k(x_1,x_2,\ldots,x_s,y)=0, \quad k=1,2,\ldots,s
$$
has a nontrivial solution in $F^{s+1}$, say
$(\xi_1,\xi_2,\ldots,\xi_s,\eta)$. It follows that the matrix
\[ \la^2 \eta^2 I_n-S(\xi_1,\ldots,\xi_s,\eta)^2 \]
is nilpotent.
We have $\eta\ne0$ since otherwise this matrix would be diagonal
and nilpotent, i.e., zero.
The matrix
\[ A=\eta^{-1}S(\xi_1,\ldots,\xi_s,\eta) \]
is bidiagonal and the matrix
$\la^2 I_n-A^2$ nilpotent. As $A$ is skew-symmetric,
its eigenvalues must be $\la$ and $-\la$, each with
multiplicity $s$.

Next we claim that $A$ is special, i.e., that all $\xi_k\ne0$.
We shall prove this claim by contradiction. So, assume that $\xi_k=0$
for some $k$. Then the matrix
$\la^2\eta^2 I_s-S_1(\xi_1,\cdots,\xi_s,\eta)$
breaks into direct sum of two blocks: the first $X_1$ of size $k$
and the second of size $s-k$ (if $k=s$ the second block is of size 0).
Similarly, the matrix
$\la^2\eta^2 I_s-S_2(\xi_1,\cdots,\xi_s,\eta)$
breaks into direct sum of two blocks: the first $X_2$ of size $k-1$
and the second of size $s-k+1$ (if $k=1$ the first block is of size 0).
Since all these four blocks are nilpotent, their traces
must be 0. Since
\begin{eqnarray*}
\tr(X_1) &=& k\la^2\eta^2+(k-1)\eta^2+(\xi_1^2+\cdots+\xi_{k-1}^2), \\
\tr(X_2) &=& (k-1)\la^2\eta^2+(k-1)\eta^2+(\xi_1^2+\cdots+\xi_{k-1}^2),
\end{eqnarray*}
we deduce that $\la^2\eta^2=0$, which is a contradiction.

Hence the matrix $\la I_n-A$ has rank $n-1$. Consequently
the elementary divisors of $A$ are indeed
$(t-\la)^s$ and $(t+\la)^s$.

The case (b) will be handled by a similar argument.
Let $S=S(x_1,x_2,\ldots,x_s,y)$ be the skew-symmetric bidiagonal
matrix whose entries on the first superdiagonal
are the indeterminates
\[ x_1,y,x_2,y,\ldots,x_{s-1},y,x_s,y \]
in this order.
The characteristic polynomial of this matrix has the form
\[ p(t)=t^{2s+1}+c_1t^{2s-1}+\cdots+c_{s-1}t^3+c_s t, \]
where
$c_k=c_k(x_1,\ldots,x_s,y)$
is a homogeneous polynomial of degree $2k$ in the indicated variables.

Since there are $s+1$ indeterminates
and $F$ is algebraically closed, the system of homogeneous equations:
$$
c_k(x_1,x_2,\ldots,x_s,y)=0, \quad k=1,2,\ldots,s
$$
has a nontrivial solution in $F^{s+1}$, say
$(\xi_1,\xi_2,\ldots,\xi_s,\eta)$. It follows that the matrix
$S(\xi_1,\ldots,\xi_s,\eta)$ is nilpotent.
We claim that $\eta\ne0$. Otherwise this matrix would be nilpotent
and semisimple, i.e., the zero matrix.

Moreover we claim that each $\xi_k\ne0$. We prove this by
contradiction. Thus assume that some $\xi_k=0$. Then the above
matrix breaks up into direct sum of two blocks: the first of
size $2k-1$ and the second of size $2(s-k+1)$.
Both of these blocks must be nilpotent. However, the determinant of
the second block is $\eta^{2(s-k+1)}\ne0$, which gives a
contradiction.

Hence the matrix
$A=\eta^{-1} S(\xi_1,\ldots,\xi_s,\eta)$
is a special bidiagonal matrix having only one elementary divisor,
$t^{2s+1}$.

The case (c) is handled by the above lemma.

\end{proof}

Next we show that there are only finitely many special
bidiagonal matrices in $\Sk_n(F)$ having prescribed
elementary divisors.

\begin{theorem} \label{spec}
If $n=2s$ is even and $\la\ne0$, then there are at most
$2^s s!$ special bidiagonal matrices $A\in\Sk_n(F)$ with
elementary divisors $(t-\la)^s$ and $(t+\la)^s$.
If $n=2s+1$ is odd, then there are at most
$2^s s!$ special bidiagonal matrices $A\in\Sk_n(F)$ with
the elementary divisor $t^{2s+1}$.
\end{theorem}

\begin{proof}
The existence of such matrices was proved in the previous theorem.
We just have to show that the system of $s$ homogeneous polynomial
equations
$c_k=c_k(x_1,\ldots,x_s,y)=0$
in $s+1$ variables has at most $2^s s!$ solutions in the associated
projective space.

We claim that the number of solutions is finite.
Otherwise the projective variety
defined by this system of equations would possess an irreducible
component, say $X$, of positive dimension. Consequently, the
intersection of $X$ with the hyperplane $y=0$ would be non-empty.
On the other hand we have shown in the proof of the previous theorem
that there are no nontrivial solutions with $y=0$.
Thus our claim is proved.

Now the assertion of the theorem follows from B\'{e}zout's Theorem
(see \cite[Chapter IV, \S2]{IS}).
\end{proof}

The non-uniqueness of special bidiagonal matrices in $\Sk_n(F)$,
with specified elementary divisors,
prevents us from obtaining a genuine normal form.

\section{Normal forms in characteristic zero}
\label{NF}

{}From now on we restrict our attention to algebraically closed
fields $F$ of characteristic 0. In this case we find very simple
normal forms for orthogonal similarity classes of
skew-symmetric matrices. In the case of nonsingular matrices,
this normal form is bidiagonal.

Let $n=2s+1$ be odd and let $R_n\in\Sk_n(F)$ be the bidiagonal
matrix whose consecutive superdiagonal entries are
$$
\sqrt{s},\,i,\sqrt{s-1},\,i\sqrt{2},\,\sqrt{s-2},\, i \sqrt 3,\,\ldots, \, \sqrt{3}, \, i\sqrt{s-2},\,
\sqrt{2},\,i\sqrt{s-1},\,1,\,i\sqrt{s}.
$$
{}For instance,
\[
R_7=\left[\begin{array}{ccccccc}
0 & \sqrt{3} & 0 & 0 & 0 & 0 & 0 \\
-\sqrt{3} & 0 & i & 0 & 0 & 0 & 0 \\
0 & -i & 0 & \sqrt{2} & 0 & 0 & 0 \\
0 & 0 & -\sqrt{2} & 0 & i\sqrt{2} & 0 & 0 \\
0 & 0 & 0 & -i\sqrt{2} & 0 & 1 & 0 \\
0 & 0 & 0 & 0 & -1 & 0 & i\sqrt{3} \\
0 & 0 & 0 & 0 & 0 & -i\sqrt{3} & 0 \end{array}\right].
\]

It is not hard to show that the matrix $R_n$ is nilpotent. The
proof is similar to the proof of \cite[Proposition 3.3]{DZ}
and we shall omit it.

{}For even $n=2s$,
let $P_n\in\Sk_n(F)$ be the bidiagonal matrix whose
consecutive superdiagonal entries are
\[ \al_k=\sqrt{\frac{k(n-k)}{(n-2k)^2-1}},\quad 1\le k\le n-1. \]
Note that for $k\ne s$ the number $\al_k^2$ is a positive rational
number while $\al_s^2=-s^2$. We may assume that $\al_s=si$.
{}For instance,
\[
P_6=\frac{1}{3} \left[\begin{array}{cccccc}
0 & \sqrt{3} & 0 & 0 & 0 & 0 \\
-\sqrt{3} & 0 & 2\sqrt{6} & 0 & 0 & 0 \\
0 & -2\sqrt{6} & 0 & 9i & 0 & 0 \\
0 & 0 & -9i & 0 & 2\sqrt{6} & 0 \\
0 & 0 & 0 & -2\sqrt{6} & 0 & \sqrt{3} \\
0 & 0 & 0 & 0 & -\sqrt{3} & 0
\end{array}\right].
\]

\begin{lemma} \label{druga}
Let $n=2s$ be even and let $P_n$ be the matrix defined above.
Then $(P_n^2-I_n)^s=0$.
\end{lemma}

The proof of this lemma is somewhat long and complicated and
will be given in the next section.

Since one has to extract square roots, there is a built in
non-uniqueness in the definition of the matrices $P_n$ and $R_n$.
Hence they should be viewed as defined only up to the choice of
these square roots, or equivalently, up to the action of
the group of diagonal matrices with the diagonal entries $\pm1$.
By abusing the language, we shall refer to this type of
non-uniqueness as the {\em choice of signs}.

Since all entries on the superdiagonal of $P_n$ are nonzero, we
conclude that the elementary divisors of $P_n$ are $(t-1)^s$
and $(t+1)^s$.

As a consequence, we obtain the following theorem.

\begin{theorem} \label{canon}
Let $F$ be an algebraically closed field of characteristic $0$.
Then any skew-symmetric matrix $A\in M_n(F)$ is orthogonally similar to
the direct sum of blocks of the following types:

(a)  $\la P_m$, $m$ even, $\la\ne0$;

(b)  $Q_m$, $m$ divisible by 4;

(c) $R_m$, $m$ odd.

This direct decomposition is unique up to the ordering of
the diagonal blocks and the choice of signs inside the blocks
of type $P_m$ and $R_m$.
\end{theorem}

We can now derive some interesting combinatorial identities
from Lemma \ref{druga}. We fix a positive integer $s$, set $n=2s$,
and define the coefficients
\begin{equation} \label{beta}
\be_k=\al_k^2=\frac{k(n-k)}{(n-2k)^2-1},\quad k\in\bZ.
\end{equation}
Note that then $\beta_k=\beta_{n-k}$ is valid for all integers $k$.

\begin{corollary} \label{ident}
For $1\le k\le s$ we have
$$
\sum_{1\le i_1\ll i_2\ll\cdots\ll i_{k-1}\ll i_k\le n-1}
\beta_{i_1}\beta_{i_2}\cdots \beta_{i_{k-1}}\beta_{i_k}=(-1)^k
{\binom{s}{k}},
$$
where $i\ll j$ means that $j-i\ge2$.
\end{corollary}

\begin{proof}
The characteristic polynomial of the matrix $P_n$ is
$$
f(t)=t^n-c_1t^{n-2}+c_2t^{n-4}-\cdots+(-1)^sc_s,
$$
where
\[
c_k=\sum_{1\le i_1\ll i_2\ll\cdots\ll i_{k-1}\ll i_k\le n-1}
\beta_{i_1}\beta_{i_2}\cdots \beta_{i_{k-1}}\beta_{i_k}.
\]
By Lemma \ref{druga}, $P_n$ has elementary divisors $(t-1)^s$, $(t+1)^s$.
Hence, $f(t)=(t^2-1)^s$ and the identities stated in the corollary follow.
\end{proof}

We end this section by proposing three related open problems.
\vskip 2mm
\noindent{\bf Problem 1}\quad Find a direct proof of the combinatorial
identities stated in the above corollary.
\vskip 2mm

\noindent{\bf Problem 2}\quad If $F$ has characteristic zero and $n=2s$ is even,
then there are exactly $2^s s!$ special bidiagonal matrices in
$\Sk_n(F)$ having elementary divisors $(t-1)^s$, $(t+1)^s$.
(We have verified this claim for $s\le4$.)
\vskip 2mm

\noindent{\bf Problem 3}\quad If $F$ has characteristic zero and $n=2s+1$
is odd, then there are exactly $2^s s!$ nilpotent special bidiagonal
matrices in $\Sk_n(F)$.
(We have verified this claim for $s\le4$.)

\section{Proof of Lemma \ref{druga}}

We recall that $n=2s$ is even and that the $\beta_k$'s
are defined by the formula (\ref{beta}).
Clearly, the matrix $P_n$ is similar to

\[
X=\left[\begin{array}{cccccc}
0 & \be_1 & 0 & & 0 & 0 \\
-1 & 0 & \be_2 & & 0 & 0 \\
0 & -1 & 0 & & 0 & 0 \\
 & & & \ddots & & \\
0 & 0 & 0 & & 0 & \be_{n-1} \\
0 & 0 & 0 & & -1 & 0
\end{array}\right].
\]

The matrix $X^2$ has zero entries in positions $(i,j)$ with
$i+j$ odd. Hence $X^2$ is permutationally similar to the
direct sum of two $s\times s$ matrices $Y$ and $Z$.
The matrix $Y$ (resp. $Z$) is the submatrix of $X^2$ occupying
the entries in positions $(i,j)$ with $i$ and $j$ odd (resp. even).
The matrices $Y$ and $Z$ are similar since $Y=UV$ and $Z=VU$,
where

\[
U=\left[\begin{array}{cccccc}
\be_1 & 0 & 0 & & 0 & 0 \\
-1 & \be_3 & 0 & & 0 & 0 \\
0 & -1 & \be_5 & & 0 & 0 \\
 & & & \ddots & & \\
0 & 0 & 0 & & \be_{n-3} & 0 \\
0 & 0 & 0 & & -1 & \be_{n-1}
\end{array}\right], \]
and
\[
V=\left[\begin{array}{cccccc}
-1 & \be_2 & 0 & & 0 & 0 \\
0 & -1 & \be_4 & & 0 & 0 \\
0 & 0 & -1 & & 0 & 0 \\
 & & & \ddots & & \\
0 & 0 & 0 & & -1 & \be_{n-2} \\
0 & 0 & 0 & & 0 & -1
\end{array}\right].
\]

It remains to prove that the matrix $Y$ is unipotent, i.e.,
the matrix $I_s-Y$ is nilpotent. For the reader's convenience,
let us display the matrix $I_s-Y$:

\[
\left[\begin{array}{cccccc}
1+\be_1 & -\be_1\be_2 & 0 & & 0 & 0 \\
-1 & 1+\be_2+\be_3 & -\be_3\be_4 & & 0 & 0 \\
0 & -1 & 1+\be_4+\be_5 & & 0 & 0 \\
 & & & \ddots & & \\
0 & 0 & 0 & & 1+\be_{n-4}+\be_{n-3} & -\be_{n-3}\be_{n-2} \\
0 & 0 & 0 & & -1 & 1+\be_{n-2}+\be_{n-1}
\end{array}\right].
\]

The proof of this last fact is quite intricate, it uses Givental's
proof of nilpotency of some special tridiagonal matrices
constructed from a finite quiver (see his paper \cite{Giv:MSFlag}).
{}For the reader's convenience, we include more detailed proof.

Denote by $\Gamma$ the infinite square grid with vertex set $\bZ^2$
and orient each horizontal edge to the right, $(i,j)\to(i+1,j)$,
and each vertical edge downward, $(i,j+1)\to(i,j)$. To this
horizontal resp. vertical edge we assign the weight
\[ u_{i,j}=-\frac{2i(2i+1)}{(2i-2j+1)(2i-2j+3)} \]
resp.
\[ v_{i,j}=\frac{2j(2j-1)}{(2i-2j-1)(2i-2j+1)}. \]
It is easy to verify that
\begin{equation} \label{Giv}
u_{i,j}+v_{i,j-1}=u_{i-1,j}+v_{i,j}, \quad
u_{i,j}v_{i,j}=u_{i,j+1}v_{i+1,j},
\end{equation}
i.e., for each vertex $(i,j)\in\Gamma$ the sum of the weights of
the two incoming edges is the same as for the two outgoing edges
and, for each small square of $\Gamma$,
the product of the edge weights along the two
oriented paths of length 2 are equal (see Figure 1).
Note that $u_{0,j}=0$ for all $j$'s and $v_{i,0}=0$ for all $i$'s.

\vspace*{2.6in}

\begin{picture}(0,0)

\put(120,0){Fig.\ 1: A portion of $\Gamma$}

\put(170,60){\line(0,1){100}}
\put(170,60){\circle*{4}}

\put(114,30){$i-1$}
\put(168,30){$i$}
\put(214,30){$i+1$}

\put(166.75,80){$\vee$}
\put(177,82){$v_{i,j-1}$}

\put(170,110){\circle*{4}}

\put(166.75,130){$\vee$}
\put(177,132){$v_{i,j}$}

\put(170,160){\circle*{4}}

\put(120,110){\line(1,0){100}}

\put(75,58){$j-1$}
\put(80,108){$j$}
\put(75,158){$j+1$}

\put(120,110){\circle*{4}}

\put(140,107.75){$>$}
\put(131,120){$u_{i-1,j}$}

\put(190,107.75){$>$}
\put(185,99){$u_{i,j}$}

\put(220,160){\circle*{4}}

\put(216.75,130){$\vee$}
\put(227,132){$v_{i+1,j}$}

\put(190,157.75){$>$}
\put(181,170){$u_{i,j+1}$}

\put(220,110){\circle*{4}}

\put(220,110){\line(0,1){50}}
\put(170,160){\line(1,0){50}}

\end{picture}

\bigskip

{}For each integer $d\ge1$ define two square matrices of size $d+1$

\[
U_d=\left[\begin{array}{cccccc}
u_{d,1} & 0 & 0 & & 0 & 0 \\
1 & u_{d-1,2} & 0 & & 0 & 0 \\
0 & 1 & u_{d-2,3} & & 0 & 0 \\
 & & & \ddots & & \\
0 & 0 & 0 & & u_{1,d} & 0 \\
0 & 0 & 0 & & 1 & 0
\end{array}\right], \]
and
\[
V_d=\left[\begin{array}{cccccc}
-1 & v_{d,1} & 0 & & 0 & 0 \\
0 & -1 & v_{d-1,2} & & 0 & 0 \\
0 & 0 & -1 & & 0 & 0 \\
 & & & \ddots & & \\
0 & 0 & 0 & & -1 & v_{1,d} \\
0 & 0 & 0 & & 0 & -1
\end{array}\right].
\]

The matrices
\[
U_d V_d=\left[\begin{array}{cccccc}
-u_{d,1} & u_{d,1}v_{d,1} & 0 & & 0 & 0 \\
-1 & v_{d,1}-u_{d-1,2} & u_{d-1,2}v_{d-1,2} & & 0 & 0 \\
0 & -1 & v_{d-1,2}-u_{d-2,3} & & 0 & 0 \\
 & & & \ddots & & \\
0 & 0 & 0 & & v_{2,d-1}-u_{1,d} & u_{1,d}v_{1,d} \\
0 & 0 & 0 & & -1 & v_{1,d}
\end{array}\right]
\]
and
\[
V_d U_d=\left[\begin{array}{cccccc}
v_{d,1}-u_{d,1} & u_{d-1,2}v_{d,1} & 0 & & 0 & 0 \\
-1 & v_{d-1,2}-u_{d-1,2} & u_{d-2,3}v_{d-1,2} & & 0 & 0 \\
0 & -1 & v_{d-2,3}-u_{d-2,3} & & 0 & 0 \\
 & & & \ddots & & \\
0 & 0 & 0 & & v_{1,d}-u_{1,d} & 0 \\
0 & 0 & 0 & & -1 & 0
\end{array}\right]
\]
are tridiagonal.

One uses induction on $d\ge1$ to prove that $U_dV_d$ is nilpotent.
The case $d=1$ is trivial to verify. Now let $d>1$.
As $U_dV_d$ and $V_dU_d$ are similar, it suffices to prove
that $V_dU_d$ is nilpotent. This is the case iff the
submatrix obtained from $V_dU_d$ by deleting the last row and
column is nilpotent.
The equalities (\ref{Giv}) imply that this submatrix is equal
to $U_{d-1}V_{d-1}$. Hence we have shown that all matrices
$U_dV_d$ are nilpotent.

{}Finally it remains to observe that $I_s-Y=U_{s-1}V_{s-1}$, which
can be easily verified by using the definition of the weights
$u_{i,j}$ and $v_{i,j}$.

\section{Some related varieties of matrices}

Let $\pN_n$ be the nilpotent cone in $M_n(F)$, i.e., the variety of
all nilpotent matrices in this algebra.
Denote by $\pS_n$ the subspace of
$M_n(F)$ consisting of all skew-symmetric bidiagonal matrices,
by $\pT_n$ the affine subspace of
$M_n(F)$ consisting of all tridiagonal matrices
$A=[a_{ij}]$ with $a_{i+1,i}=1$ for $i=1,\ldots,n-1$,
and by $\pB_n$ the affine subspace of $\pT_n$
consisting of the matrices having zero diagonal.
Finally, let
\[ \pS^*_n\subseteq\pS_n,\quad
\pT^*_n\subseteq\pT_n,\quad
\pB^*_n\subseteq\pB_n \]
be the open subvarieties consisting of the matrices having nonzero
entries along the first super-diagonal.
Note that any matrix in $\pS_n^*$ is similar to a unique matrix
in $\pB_n^*$, and the
resulting map $\pS_n^*\to\pB_n^*$ is a $2^{n-1}$--fold covering.

Unless stated otherwise, we assume from now on that $n=2s+1$ is odd.
We are interested in the intersections 
$\pB_n\cap\pN_n$ and $\pS_n\cap \pN_n$ and their coordinate rings. 
We will relate these to a closed subvariety $\pV_s\subseteq F^{n-1}$
to be defined shortly.
We introduce first three maps
\[ A_1:F^{n-1}\to\pT_{s+1}, \quad
A_2:F^{n-1}\to\pT_s, \quad
B:F^{n-1}\to\pB_n, \]
by defining
\[ A_1(\bp)=\left[
\begin{array}{ccccc}
 -p_1 & p_1 p_2 & & & \\
 1 & -p_2-p_3 & p_3 p_4 & & \\
  & 1 & -p_4-p_5 & & \\
  & & & \ddots & p_{n-2}p_{n-1} \\
  & & & 1 & -p_{n-1}
\end{array} \right], \]
\[ A_2(\bp)=\left[
\begin{array}{ccccc}
 -p_1-p_2 & p_2 p_3 & & & \\
 1 & -p_3-p_4 & p_4 p_5 & & \\
  & 1 & -p_5-p_6 & & \\
  & & & \ddots & p_{n-3}p_{n-2} \\
  & & & 1 & -p_{n-2}-p_{n-1}
\end{array} \right], \]
and
\[ B(\bp)=\left[
\begin{array}{ccccc}
 0 & -p_1 & & & \\
 1 & 0 & -p_2 & & \\
 & 1 & 0 & & \\
 & & & \ddots & -p_{n-1} \\
 & & & 1 & 0
\end{array} \right], \]
where $\bp=(p_1,\ldots,p_{n-1})\in F^{n-1}$.

\begin{proposition} \label{p:iso}
The characteristic polynomials
of $A_1(\bp)$, $A_2(\bp)$ and $B(\bp)$ are related as follows.
\begin{eqnarray}\label{e:detA1}
\det(\, t\; I_{s+1}- A_1(\bp))&= &t\, \det (\, t\; I_{s}- A_2(\bp)),\\
\det(\, t\; I_{n}- B(\bp))&=  & t\, \det (\, t^2\; I_{s}- A_2(\bp)).
\label{charpolB}
\end{eqnarray}
Hence, if one of the matrices
$A_1(\bp)$, $A_2(\bp)$, $B(\bp)$ is nilpotent, so are
the other two. 
\end{proposition}
\begin{proof}
Following Givental \cite{Giv:MSFlag}, we define the matrices (of size $s+1$)
\[ U(\bp)=\left[
\begin{array}{cccccc}
 -p_1 & & & & & \\
 1 & -p_3 & & & & \\
 & 1 & -p_5 & & & \\
 & & & \ddots & & \\
 & & & & -p_{n-2} & \\
 & & & & 1 & 0
\end{array} \right] \]
and
\[ V(\bp)=\left[
\begin{array}{cccccc}
 1 & -p_2 & & & & \\
   & 1 & -p_4 & & & \\
   & & 1 & & & \\
 & & & \ddots & & \\
 & & & & 1 & -p_{n-1} \\
 & & & & & 1
\end{array} \right]. \]
Then $U(\bp)V(\bp)=A_1(\bp)$ and
\[ V(\bp)U(\bp)=
\left[ \begin{array}{cc} A_2(\bp) & 0 \\ w & 0 \end{array} \right], \]
where $w$ is a row vector. Consequently,
\begin{multline*}
\det(\, t\; I_{s+1}- A_1(\bp))=
\det(\, t\; I_{s+1}- U(\bp) V(\bp))\\
=\det(\, t\; I_{s+1}-V(\bp) U(\bp))=
t \det(\, t\; I_{s}-A_2(\bp)).
\end{multline*}

It is easy to verify that there is a permutation matrix
$\Pi_n\in M_n(F)$ such that
\[ \Pi_n B(\bp)^2 \Pi_n^{-1} =
\left[ \begin{array}{cc} A_1(\bp)  & 0 \\ 0 & A_2(\bp)
\end{array} \right] \]
holds for all $\bp\in F^{n-1}$.
So, by using \eqref{e:detA1},
\begin{multline*}
\det(\, t\; I_{n}- B(\bp)^2)=
\det(\, t\; I_{s+1}- A_1(\bp))
\det(\, t\; I_{s}-A_2(\bp))\\
= t \det(\, t\; I_{s}-A_2(\bp))^2.
\end{multline*}
On the other hand it is easy to see that
$
\det(\, t\; I_{n}- B(\bp))=
\det(\, t\; I_{n}+B(\bp))
$
and so
$$
\det(\, t\; I_{n}- B(\bp))^2
=\det(\, t^2\; I_{n}- B(\bp)^2)=
t^2 \det(\, t^2\; I_{s}-A_2(\bp))^2.
$$
Therefore $\det(\, t\; I_{n}- B(\bp))=\pm
t \det(\, t^2\; I_{s}-A_2(\bp))$ which
after comparing the leading coefficients implies
\eqref{charpolB}.
\end{proof}

The non-constant coefficients of the characteristic polynomial
of $B(\bp)$, or equivalently of $A_1(\bp)$ or $A_2(\bp)$ define a
 variety (as we will see, reduced)
inside $F^{2s}$ which we denote by $\mathcal V_{s}$. In other words,
\[ \pV_s=A_1^{-1}(\pT_{s+1}\cap\pN_{s+1})=A_2^{-1}(\pT_s\cap\pN_s)=
B^{-1}(\pB_n\cap\pN_n). \]
Define also $\pV_s^*=\pV_s\cap(F^*)^{n-1}$ and note that
the map $B$ induces isomorphisms
$\pV_s\to\pB_n\cap\pN_n$ and $\pV_s^*\to\pB_n^*\cap\pN_n$.

Let $\al_1:\pV_s\to\pT_{s+1}\cap\pN_{s+1}$
and $\al_1^*:\pV_s^*\to\pT^*_{s+1}\cap\pN_{s+1}$
be the maps induced by $A_1$,
and define similarly the maps $\al_2$ and $\al_2^*$.

\begin{proposition} \label{p:components}
$\pV_{s}^*$ is a smooth, irreducible variety of dimension
$s$ and $\al_1^*$ is an open inclusion.
\end{proposition}

\begin{proof}
By Kostant's work on the Toda lattice (see \cite[Theorem 2.5]{Kos:Toda})
the scheme-theoretic intersection $\pT^*_{s+1}\cap\pN_{s+1}$ 
defines a smooth, irreducible variety of
dimension $s$. In fact it is isomorphic to an open subset of $F^s$.
It suffices, therefore, to show that $\al_1^*$ is an open embedding. A matrix
\[ \left[ \begin{array}{ccccc}
-a_1 & b_1 & & & \\
 1 & -a_2 & b_2 & & \\
 & 1 & \ddots & & \\
 & & & -a_s & b_s \\
 & & & 1 & -a_{s+1}
\end{array} \right]
\in \pT^*_{s+1}\cap\pN_{s+1} \]
lies in the image of $\al_1^*$ precisely if the denominators in
the continued fraction expansions
\begin{align*}
p_1&=a_1,  & & p_2&=\frac{b_1}{a_1} ,& &
p_3&=a_2+\frac{b_1}{a_1} ,& &
p_4&=\frac{b_2}{a_2+\frac{b_1}{a_1}}, & &
p_5&=a_3+\frac{b_2}{a_2+\frac{b_1}{a_1}}
\end{align*}
up to
\[
p_{n-1}=\frac{b_s}{a_s+\frac{b_{s-1}}{\cdots +\frac{b_1}{a_1}}}
\]
are all nonzero. This clearly defines an open subset of
$\pT^*_{s+1}\cap\pN_{s+1}$. The above formulas for the $p_i$
define an algebraic  inverse from this open set to $\pV_s^*$.
\end{proof}

We have therefore proved that $\pB^*_n\cap\pN_n$
is a smooth, irreducible variety of dimension $s$.

Denote by $\Lambda_s$ the set of subsequences
$\la=(\la_1,\ldots,\la_k)$ of $(1,2,\ldots,2s)$ with $k$ even,
$0\le k\le 2s$, and such that $\la_i\equiv i \pmod{2}$ for each $i$.
Let $\bp=(p_1,\ldots,p_{2s})\in\pV_s$ and let $\la_\bp=(\la_1,\ldots,\la_k)$
be the increasing sequence of indices such that
\[ p_{\la_1}=\cdots=p_{\la_k}=0 \]
and all other coordinates $p_i$ are nonzero. Since nilpotent matrices have
zero determinant, it is easy to see that $\la_\bp\in\Lambda_s$.
Consequently, we have a map $\sigma:\pV_s\to\Lambda_s$ defined by
$\sigma(\bp)=\la_\bp$. We denote by $\pV_s^\la$ the fibre of $\sigma$
lying over the point $\la\in\Lambda_s$. This gives a set-theoretic
partition
\begin{equation} \label{part}
    \pV_s=\coprod_{\la\in\Lambda_s} \pV_s^\la.
\end{equation}
These fibres are in fact smooth subvarieties. For instance, for the
empty sequence $\emptyset\in\Lambda_s$ we have $\pV_s^\emptyset=\pV_s^*$.
Moreover, we have the following description of the fibres.

\begin{proposition} \label{fibres}
For $\la=(\la_1,\ldots,\la_k)\in\Lambda_s$,
$k=2r$, we set $\la_0=0$, $\la_{k+1}=n$ and  $s_i=(\la_i-\la_{i-1}-1)/2$,
$1\le i\le k+1$. Then
\[ \pV_s^\la\cong\prod_{1\le i\le k+1} \pV^*_{s_i}, \]
where, by convention, $\pV_0^*$ is the variety consisting of a single point.
In particular $\pV_s^\la$ is a smooth variety of dimension $s-r$.
\end{proposition}

\begin{proof}
Let $\bp\in \pV_s^\la$. Then the characteristic polynomial
of $B(\bp )$ agrees with the characteristic polynomial of the block matrix
\begin{equation*}
M=
\left[
\begin{array}{ccccc}
M_1 & & & &\\
& M_2& & &\\
& & & \ddots&\\
& & & &M_{k+1},
\end{array}
\right ],
\end{equation*}
where
\begin{equation*}
M_i=
\left[
\begin{array}{ccccc}
 0 & -p_{\la_{i-1}+1} & & & \\
 1 & 0 & -p_{\la_{i-1}+2} & & \\
 & 1 & 0 & & \\
 & & & \ddots & -p_{\la_i-1} \\
 & & & 1 & 0
\end{array} \right].
\end{equation*}
The matrix $M$ is nilpotent if and only if each of the $M_i$ is
nilpotent. Therefore
\[\bp_i:=(p_{\la_{i-1}+1},p_{\la_{i-1}+2},\dotsc, p_{\la_{i}-1})\]
lies in $\pV^*_{s_i}$. The resulting map
\begin{eqnarray*}
\pV_s^\la &\to &\prod_{1\le i\le k+1} \pV^*_{s_i} \\
 \bp& \mapsto & (\bp_1,\bp_2,\dotsc, \bp_{k+1})
\end{eqnarray*}
is clearly an isomorphism.
\end{proof}

\begin{theorem} \label{irr}
$\pV_s$ is an irreducible 
%normal 
variety.
\end{theorem}

\begin{proof}
$\pV_s$ is defined
as the intersection
of $s$ hypersurfaces, the zero sets of the coefficients of $t^{k}$,
$k=0,\ldots,s-1$, of the characteristic polynomial of $A_2(\bp)$.
Consequently, each irreducible component of $\pV_s$ has dimension at
least $s$. On the other hand, Proposition~\ref{fibres} implies that
these dimensions are at most equal to $s$. Hence, $\pV_s$ is an
equidimensional variety of dimension $s$. It follows that
$\pV_s$ is irreducible.
\end{proof}

\begin{corollary}\label{skewvar}
$\pS_n\cap \pN_n$
and $\pS^*_n\cap\pN_n$ are varieties of pure dimension $s$.
Moreover
$\pS^*_n\cap\pN_n$ is smooth and an open dense subset 
of $\pS_n\cap\pN_n$.
\end{corollary}

\begin{proof}
We have a well-defined map $\phi:\pS_n\cap\pN_n\to \pV_s$
which takes a nilpotent skew-symmetric bidiagonal matrix with entries
$(a_1,\dotsc, a_{n-1})$ along the upper diagonal, to an element
$(-a_1^2,\dotsc, -a_{n-1}^2)\in \mathcal V_s$.  This map restricts to 
a covering over each of the fibers $\mathcal V_s^{\lambda}$, which implies that 
$\pS_n\cap\pN_n$ is reduced and any irreducible component 
has dimension
at most $s$. However, as for $\mathcal V_s$, any irreducible component 
of $\pS_n\cap\pN_n$ has
dimension at least $s$. Therefore $\pS_n\cap\pN_n$ is of pure dimension $s$. 
As a consequence the $s$-dimensional fiber, $\pS^*_n\cap\pN_n$, is an open dense subset
of $\pS_n\cap\pN_n$. Moreover, since $\pS^*_n\cap\pN_n$ is a covering space of 
$\mathcal V_s^*$, it is a smooth variety. 
\end{proof}

The problem of finding ``nice'' bidiagonal
normal forms for orthogonal similarity classes of skew-symmetric
matrices is identical to exhibiting ``nice'' representatives of
the intersection $\pS_n\cap\pN_n$.
In connection with this it is important to determine
precisely the  irreducible components of $\pS_n\cap\pN_n$.
If $n=3$ it is easy to check directly that $\pS_3\cap\pN_3$
has two irreducible components. 

\vskip 2mm

\noindent{\bf Problem 4}\quad Prove or disprove the following assertion.
For $n>3$ odd, the variety $\pS^*_n\cap\pN_n$ is connected,
 and $\pS_n\cap\pN_n$ is irreducible.

\begin{remark}\label{r:evencase} The above considerations 
are for varieties of $n\times n$ matrices where $n$ is odd.  
For even $n=2s$ we remark that  
$\mathcal B_{2s}\cap\mathcal N_{2s}$ 
is a variety of pure dimension $s-1$, but no longer irreducible if $s>1$. 
Namely it is easy to see that one has a  decomposition into $s$ 
components, 
$$
\mathcal B_{2s}\cap\mathcal N_{2s}=\bigcup_{j=1}^s
\mathcal I_j,
$$
where $\mathcal I_j$ is defined by the vanishing of the 
$(2j-1)^{st}$ entry on the upper diagonal. Each of these components is 
isomorphic to a product, 
$
\mathcal I_j\cong \mathcal V_{j-1}\times \mathcal V_{s-j},
$
and hence is irreducible and $(s-1)$-dimensional.

For the even-dimensional skew-symmetric case the scheme-theoretic intersection 
$\pS_{2s}\cap\pN_{2s}$ is not reduced, so we choose to consider
$\pS_{2s}\cap\pN_{2s}$ as an intersection of algebraic sets. The resulting variety 
then  
again decomposes into a union of $s$ subvarieties $\mathcal I_j^{\rm skew}$, 
where the 
definition of $\mathcal I_j^{\rm skew}$ is analogous to that of $\mathcal I_j$ above. 
Moreover, the component $\mathcal I_j^{\rm skew}$ is again isomorphic to a product,
\begin{equation}\label{e:subdivision}
(\pS_{2j-1}\cap\pN_{2j-1}) \times (\pS_{2s-2j+1}\cap\pN_{2s-2j+1}),
\end{equation}
and Corollary~\ref{skewvar} implies that it has pure dimension $s-1$. Hence 
$\pS_{2s}\cap\pN_{2s}$ has pure dimension $s-1$. 

The number of irreducible components of $\pS_{2s}\cap\pN_{2s}$, when $s>1$,
depends on the numbers of irreducible components of the subvarieties
\eqref{e:subdivision}. 
Namely, since $\pS_{3}\cap\pN_{3}$ has two irreducible
components and thus 
$\mathcal I_2^{\rm skew}$ and $\mathcal I_{s-1}^{\rm skew}$ 
have at least two irreducible components each, the variety 
$\pS_{2s}\cap\pN_{2s}$ must have at least $s+2$ irreducible components. 

Explicitly we have the following. In the case $s=1$ the algebraic set 
$\pS_{2}\cap\pN_{2}$ consists of a single 
point. For $s=2$ there are precisely 
$s+2=4$ irreducible components in
$\pS_{4}\cap\pN_{4}$, obtained as described above. 
The case $s=3$ is exceptional with 
$6$ irreducible components. This is because the variety
$\mathcal I^{\rm skew}_2$ is isomorphic to 
$(\pS_{3}\cap\pN_{3})\times (\pS_{3}\cap\pN_{3})$ which
has $4$ irreducible components, and it can be checked that 
$\mathcal I^{\rm skew}_1$ and $\mathcal I^{\rm skew}_3$  are irreducible.   
For $s\ge 4$, if the assertion in Problem~4 is true then
the variety $\pS_{2s}\cap\pN_{2s}$ 
has precisely $s+2$ irreducible components.   
\end{remark}

The results in this section may be compared with the normalized tridiagonal
and the tridiagonal symmetric cases.

In the normalized tridiagonal case $\pT_m\cap\pN_m$ was shown to be irreducible
(of dimension $m-1$) by Kostant \cite{Kos:QCoh}. This variety is of
particular interest as its coordinate ring has another interpretation
as the quantum cohomology ring $qH^*(\SL_m/B, F)$ of the flag variety
$\SL_m/B$.
In this context, the coordinate ring of our variety $\pB_n\cap\pN_n$
can also be interpreted as the quotient of the quantum cohomology
ring $qH^*(\SL_n/B,F)$ by the ideal generated by the Chern classes of
the tautological line bundles,
$x_i=c_1(L_i)$. Thus Theorem~\ref{irr} implies that
$$
qH^*(\SL_{n}/B,F)/(x_1,\dotsc, x_n)
$$
is an integral domain (for odd $n$).

The variety of symmetric nilpotent $m\times m$ matrices is also
equidimensional, as can be shown without much difficulty,  and is
irreducible precisely if its intersection with the regular nilpotent
orbit (which is smooth) is connected.


\begin{thebibliography}{99}

\providecommand{\bysame}{\leavevmode\hbox to3em{\hrulefill}\thinspace}

\bibitem{DZ}
D.\v{Z}. \DJo\; and K. Zhao, Tridiagonal normal forms
for orthogonal similarity classes of symmetric matrices,
Linear Algebra and its Applications 384 (2004) 77--84.

\bibitem{FG}
F.R. Gantmacher, The Theory of Matrices, vol. 2, Chelsea,
New York, 1989.

\bibitem{Giv:MSFlag}
A. Givental, Stationary phase integrals, quantum {T}oda lattices, flag
manifolds and the mirror conjecture,
Topics in singularity theory, American Mathematical Society
Translations Ser 2., AMS, 1997.



\bibitem{Kos:Toda}
B. Kostant, The solution to a generalized {T}oda lattice and
representation theory, Adv. in Math. \textbf{34} (1979), no.~3, 195--338.

\bibitem{Kos:QCoh}
\bysame, Flag manifold quantum cohomology, the {T}oda lattice, and the
representation with highest weight $\rho$,
Selecta Math. (N.S.) \textbf{2} (1996), 43--91.




\bibitem{IK}
I. Kaplansky, Linear Algebra and Geometry, A second course,
Allyn and Bacon, Boston, 1969.


\bibitem{IS}
I.R. Shafarevich, Basic Algebraic Geometry, Nauka, Moscow,
1972. English transl.: Grundlehren der mathematischen
Wissenschaften 213, Springer, New York, Berlin, Heidelberg, 1974.

\end{thebibliography}
\end{document}